\documentclass[11pt,twoside,reqno]{article}
\usepackage{amsfonts}
\usepackage{graphicx}
\begin{document}
\begin{center}
{\Large{\textbf{ILU Preconditioning Based on the FAPINV Algorithm}}}
\\
\vspace{0.5cm}

{\bf Davod Khojasteh Salkuyeh}

{\small Department of Mathematics, University of Mohaghegh Ardabili,\\
P. O. Box. 179, Ardabil, Iran\\
E-mail: khojaste@uma.ac.ir}

\bigskip

{\bf Amin Rafiei}

{\small Department of Mathematics, Sabzevar Tarbiat Moallem University,\\
P. O. Box. 397 , Sabzevar, Iran\\
E-mail: rafiei.am@gmail.com, rafiei@sttu.ac.ir}

\bigskip

{\bf Hadi Roohani}

{\small Department of Mathematics, Imam Khomeini International University,\\ Ghazvin, Iran\\
E-mail: hadiroohani61@gmail.com}

\end{center}

\bigskip

\medskip

\begin{center}
 {\bf Abstract}
\end{center}

A technique for computing an ILU preconditioner based on the FAPINV algorithm is presented.
We show that this algorithm is well-defined for H-matrices.  Moreover, when used in conjunction with Krylov-subspace-based iterative solvers such as the GMRES algorithm, results in reliable solvers. Numerical experiments on some test matrices are given to show the efficiency of the new ILU preconditioner.


\bigskip

\bigskip

\bigskip

\noindent {\bf 1. Introduction}

\bigskip

Consider the linear system of equations
\begin{equation}\label{e:eq01}
    Ax=b,\\
\end{equation}
where the coefficient matrix  $A \in \mathbb{R}^{n \times n}$ is nonsingular, large, sparse and $x,b\in \mathbb{R}^{n}$. Such linear systems are often solved by Krylov subspace methods such as the
GMRES \cite{GMRES} and the BiCGSTAB \cite{BiCGSTAB} methods. In general, the convergence of the Krylov subspace methods is not guaranteed or it may be extremely slow. Hence, the original system (\ref{e:eq01}) is transformed into a more tractable form.  More precisely, to obtain good convergence rates, or even to converge, Krylov subspace methods are applied to the left preconditioned linear system
\[
MAx=Mb,
\]
\noindent or to the right preconditioned linear system
\[
AMy=b,\quad x=My,
\]
where the matrix $M$ is a proper preconditioner. Two-side preconditioning is also possible \cite{Survey,Comp,Saadbook}.

There are two general ideas for constructing a preconditioner. The first one is to find a matrix $M$ such that $M=G^{-1}$, where $G$ approximates $A$ in some sense. In this case, $M$ should be chosen such that $AM$ (or $MA$) is a good approximation of the identity matrix.
The best-known general-purpose preconditioners in this class are those based on the incomplete LU (ILU) factorization of the original matrix. Let the matrix $A$ admits the LDU factorization
$A=LDU$, where $L$ and $U^T$ are lower unitriangular matrices and  $D={\rm diag}(d_1,d_2,\ldots,d_n)$ is a diagonal matrix. Here $G=\tilde{L} \tilde{D} \tilde{U}$, where  $\tilde{L}$ and $\tilde{U}^T$ are sparse lower unitriangular matrices which approximate $\tilde{L}$ and $\tilde{U}^T$, respectively, and  $\tilde{D}$ is a diagonal matrix which approximates $D$. The ILU preconditioners are very much effective in increasing the rate of convergence. Their main drawback is the possibility of breakdowns during incomplete factorization process, due to the occurrence of zero or small pivots (or appearing nonpositive pivot elements for symmetric positive definite (SPD) matrices). In \cite{Meij}, Meijerink and van der Vorst  have shown that this type of preconditioning exists for M-matrices. Then, Manteuffel in \cite{Mant} has extended this result to the H-matrices. We recall that the matrix $A=(a_{ij})$ is an M-matrix if $a_{ij}\leq 0$ for all $i \neq j$, $A$ is nonsingular and $A^{-1} \geq 0$. Moreover, $A$ is an H-matrix if its comparison matrix $\hat{A}=(\hat{a}_{ij})$ is an M-matrix where
\[\hat{a}_{ij}=
\left\{
  \begin{array}{ll}
    -|a_{ij}|, & i\neq j, \\
    \hspace{0.3cm}|a_{ii}|, & i=j.
  \end{array}
\right.
\]
For general matrices, there are some ways to guard against appearing zero or very small pivots, see for example \cite{RIF,Rezghi,Lee1}. Another drawback of the ILU preconditioners is the lack of inherent parallelism. Many researchers have made effort to improve the accuracy and the degree of parallelism of the ILU preconditioners in the past \cite{Saadbook,ILUT,ILUM,BILUM,BILUTM,Benzi1,Benzi2,Meurant}.

The second idea for constructing a preconditioner is to find a matrix $M$ that directly approximates $A^{-1}~(M \approx A^{-1}$). In this case, in practice we do not need to compute $AM~(\textrm{or}~MA)$ explicitly, because when the Krylov subspace methods are used to solve a preconditioned system, only the matrix-vector product is required. One drawback of many sparse approximate inverse techniques is their high construction cost, unless the computation can be done efficiently on parallel computers. One approach in this class is to compute  sparse approximate inverse in the factored form. Here from LDU factorization $A=LDU$, we have $A^{-1}=ZD^{-1}W$ where $W=L^{-1}$ and $Z=U^{-1}$. If $G=\tilde{Z} \tilde{D} \tilde{W}$, where $\tilde{W}\approx W,~ \tilde{Z} \approx Z$, and $\tilde{D} \approx D^{-1}$, in which $\tilde{W}$ and $\tilde{Z}^T$ are sparse lower unitriangular matrices and $\tilde{D}$ is diagonal matrix, then $G$ may be used as a preconditioner for system (\ref{e:eq01}) and is called factored sparse approximate inverse.  Here $\tilde{W}$ and $\tilde{Z}$ are called  sparse approximate inverse factors of $A$. There are several algorithms to compute a factored sparse approximate inverse of a matrix. Among them are the FSAI algorithm \cite{Kol1,Kol2}, the AIB algorithm \cite{Saadbook}, the AINV method \cite{SYMAINV,AINV}, and the FAPINV algorithm \cite{Luo1,Luo2,Luo3,ZhangJi,Lee2,Zhang}.

For SPD matrices, there exists a variant of the AINV method, denoted by SAINV (for Stabilized AINV), that is breakdown-free \cite{SAINV}. This algorithm is also presented with the name of AINV-A, independently,  by Kharchenko et al. in \cite{AINVA}. Benzi and Tuma in \cite{RIF} have introduced an ILU factorization based on the SAINV algorithm. In the proposed algorithm the $L$ factor of $\rm {LDL^T}$ factorization of  $A$ can be obtained as a by-product of the $A$-orthogonalization process used in the SAINV algorithm, at no extra cost.  Rezghi and Hosseini in \cite{Rezghi}, have shown that a similar algorithm free from breakdown can be established for nonsymmetric positive definite matrices.

The main idea of the FAPINV algorithm was introduced by Luo \cite{Luo1,Luo2,Luo3}.
Then the algorithm was more investigated by Zhang in \cite{Zhang}. Since in this procedure the factorization is performed in backward direction, we call it BFAPINV (Backward FAPINV) algorithm. In \cite{ZhangJi}, Zhang proposed an alternative procedure to compute the factorization in the forward direction, which we call it FFAPINV (Forward FAPINV) algorithm. In \cite{Lee2}, Lee and Zhang have shown that the BFAPINV algorithm is well-defined for M-matrices. It can be easily seen that this is correct for the FFAPINV algorithm as well. In \cite{Salkuyeh}, Salkuyeh showed that the FFAPINV algorithm with a simple revision may be used  for the nonsymmetric positive definite matrices, free from breakdown.

In this paper, we show that the FFAPINV algorithm is free from breakdown for H-matrices and we propose a technique for computing an ILU preconditioner based on the FAPINV algorithm at no extra cost.

\indent This paper is organized as follows. In section 2,  we review the FFAPINV algorithm.
Section 3 is devoted to the main results. Numerical experiments are given in
section 4. Finally, we give some concluding remarks in section 5.

\bigskip

\bigskip

\medskip

\noindent  \textbf{2. A review  of the FFAPINV algorithm}

\bigskip

\medskip

\indent Let $W$ and $Z^T$ be lower unitriangular matrices and $D$ be a diagonal matrix. Also, suppose that
\begin{equation}\label{e:eq02}
WAZ=D^{-1}.
\end{equation}
\noindent In this case, we term $W,Z$ and $D$, the inverse factors of  $A=(a_{ij})$. Consider $W=(w_1^T,w_2^T,\ldots,w_n^T)^T$, $Z=(z_1,z_2,\ldots,z_n)$ and $D={\rm diag}(d_1,d_2,\ldots,d_n)$, in which $w_i$'s and $z_i$'s are the rows and columns of $W$ and $Z$, respectively. Using Eq. (\ref{e:eq02}) we obtain
\begin{equation}\label{e:eq03}
w_iAz_j=\left\{
  \begin{array}{ll}
    \frac{1}{d_i}, & i=j, \\
    0, & i\neq j.
  \end{array}
\right.
\end{equation}
\noindent From the structure of the matrices $W$ and $Z$, we have
\begin{equation}\label{e:eq04}
 z_1=e_1, \quad z_j = e_j-\sum_{i=1}^{j-1}\alpha_{i}^{(j)} z_i,  \quad j=2,\ldots,n,
\end{equation}
\begin{equation}\label{e:eq05}
  w_1=e_1^T, \quad w_j = e_j^T-\sum_{i=1}^{j-1}\beta_{i}^{(j)} w_i,  \quad j=2,\ldots,n,
\end{equation}
\noindent for some $\alpha_{i}^{(j)}$'s and $\beta_{i}^{(j)}$'s, where $e_j$ is the $j$th column of the identity matrix.

First of all, from (\ref{e:eq03}) we see that
\[
d_1=\frac{1}{z_1^TAz_1}=\frac{1}{e_1^TAe_1}=\frac{1}{a_{11}}.
\]
Now let $2\leq j \leq n$ be fixed. Then from  Eqs. (\ref{e:eq03}) and (\ref{e:eq04}) and for $k=1,\ldots,j-1$, we have
\begin{eqnarray*}
  0 \hspace{-0.2cm}&=&\hspace{-0.2cm} w_k A z_j = w_k A e_j-\sum_{i=1}^{j-1}\alpha_{i}^{(j)} w_k A z_i
    = w_k A_{*j}-\alpha_{k}^{(j)} w_k A z_k= w_k A_{*j}-\alpha_{k}^{(j)}\frac{1}{d_k},
\end{eqnarray*}
where $A_{*j}$ is the $j$th column of $A$.
Therefore
\[
\alpha_{i}^{(j)}=d_iw_iA_{*j},\quad i=1,\ldots,j-1.
\]
In the same manner
\[
\beta_{i}^{(j)}=d_i A_{j*}z_i,\quad i=1,\ldots,j-1,
\]
where $A_{j*}$ is the $j$th row of $A$.  Pre-multiplying both sides of (\ref{e:eq04}) by $w_jA$ yields
\[
d_j=\frac{1}{w_jA_{*j}}.
\]
Putting these results together gives the following algorithm (FFINV for Forward Factored INVerse) for computing the inverse factors of $A$.

\bigskip

\underline{{\bf Algorithm 1:} FFINV algorithm (vector form)}
\vspace{-0.3cm}

\begin{enumerate}
\item  $z_1:=e_1$, $w_1:=e_1^T$ and $d_1:=\frac{1}{a_{11}}$ \vspace{-0.25cm}

\item For $j=2,\ldots,n,$~ Do \vspace{-0.25cm}

\item \qquad  $z_j:=e_j$; \quad $w_j:=e_j^T$ \vspace{-0.25cm}

\item \qquad  For $i=1,\ldots,j-1,$ Do \vspace{-0.25cm}

\item \qquad\qquad  $\alpha_{i}^{(j)}:=d_iw_iA_{*j}$; \quad  $\beta_{i}^{(j)}:=d_iA_{j*}z_i$\vspace{-0.25cm}

\item \qquad\qquad  $z_j:=z_j-\alpha_{i}^{(j)} z_i$; \quad  $w_j:=w_j-\beta_{i}^{(j)} w_i$ \vspace{-0.25cm}

\item \qquad  EndDo \vspace{-0.25cm}

\item \qquad  $d_j:=\frac{1}{w_jA_{*j}}$ \vspace{-0.25cm}

\item  EndDo
\item Return $W=[w_1^T,\cdots,w_n^T]^T$, $Z=[z_1,\cdots,z_n]$ and $D=diag(d_1,\cdots,d_n)$.
\end{enumerate}

Algorithm 1, is the vector form of the FFINV algorithm. It can be easily verified that this algorithm is equivalent to the following algorithm (see \cite{Zhang}). Moreover, the values of  $\alpha_{i}^{(j)}$'s and $\beta_{i}^{(j)}$'s are the same in both Algorithms 1 and 2. In this algorithm, we assume $w_j=(w_{j1},w_{j2},\ldots,w_{jn})$ and $z_j=(z_{1j},z_{2j},\ldots,z_{nj})^T$.\\

\underline{{\bf Algorithm 2:} FFINV algorithm}
\vspace{-0.25cm}

\begin{enumerate}
\item $W:=I_{n \times n}$, $Z:=I_{n \times n}$.  \vspace{-0.25cm}

\item For $j=1,\ldots,n,$ Do \vspace{-0.25cm}

\item \qquad For $i=j-1,\ldots,1,$ Do\vspace{-0.25cm}

\item \qquad\qquad $l_{i}^{(j)}=a_{ji}+\sum^{i-1}_{k=1}a_{jk} z_{ki}$ \vspace{-0.25cm}

\item \qquad\qquad $\beta_{i}^{(j)}=l_{i}^{(j)}d_i$ \vspace{-0.25cm}

\item \qquad EndDo\vspace{-0.3cm}

\item \qquad For $i=1,\ldots,j-1,$ Do\vspace{-0.25cm}

\item \qquad\qquad $w_{ji}=-\beta_{i}^{(j)}-\sum^{j-1}_{k=i+1} \beta_{k}^{(j)} w_{ki}$\vspace{-0.25cm}

\item \qquad EndDo,\vspace{-0.25cm}

\item \qquad $d_{j}={1}/{(a_{jj}+\sum^{j-1}_{k=1}w_{jk} a_{kj})}$\vspace{-0.25cm}

\item \qquad For $i=j-1,\ldots,1,$ Do \vspace{-0.3cm}

\item \qquad\qquad $u_{i}^{(j)}=a_{ij}+\sum^{i-1}_{k=1}w_{ik}a_{kj} $ \vspace{-0.25cm}

\item \qquad\qquad $\alpha_{i}^{(j)}=u_i^{(j)}d_i$ \vspace{-0.25cm}

\item \qquad EndDo\vspace{-0.25cm}

\item \qquad For $i=1,\ldots,j-1$ Do\vspace{-0.25cm}

\item \qquad\qquad $z_{ij}=-\alpha_{i}^{(j)}-\sum^{j-1}_{k=i+1}\alpha_{k}^{(j)}z_{ik}$,\vspace{-0.25cm}

\item \qquad EndDo\vspace{-0.25cm}
\item EndDo
\item Return $W=(w_{ij})$, $Z=(z_{ij})$ and $D=diag(d_1,\cdots,d_n)$.
\end{enumerate}

This algorithm computes the inverse factors $W,Z$ and $D$ such that Eq. (\ref{e:eq02}) holds.
Therefore, we have $A^{-1}=ZDW$. This shows that the inverse of $A$ in the factored form can be computed by this algorithm. A sparse approximate inverse of $A$  is computed by inserting some dropping strategies in Algorithm 2.
A dropping strategy can be used  as follows. The dropping strategy is applied in four places in the algorithm. In step 5, if $\mid \beta_i^{(j)} \vert<\tau$, then $\beta_i^{(j)}:=0$, and in step 8, if $\vert w_{ji}\vert<\tau $, then $w_{ji}:=0$. In the same way $\alpha_i^{(j)}$ in step 13 and $z_{ij}$ in step 16 are dropped when their absolute values are less than tolerance $\tau$. Hereafter, FFAPINV algorithm refers to FFINV algorithm with this type of dropping.
Here we mention that the dropping strategy proposed in \cite{ZhangJi} is slightly different from our dropping strategy.   In the next section, we show that the FFAPINV algorithm is well-defined for H-matrices.

\bigskip

\bigskip

\medskip

\noindent  \textbf{3. Existence of FFAPINV algorithm for H-matrices}

\bigskip

\medskip

First we state the following theorem.\\

\noindent \textbf{Theorem 1.} \textit{Assume that $A$ is an  M-matrix. Let $W$ and $Z$ be the inverse factors of $A$ computed by FFINV algorithm, i.e., $WAZ=D^{-1}$. Also suppose that and $\hat{W}$ and $\hat{Z}$ be the  inverse factors of $A$  computed by the FFAPINV algorithm, i.e., $\hat{W} A \hat{Z}\approx \hat{D}^{-1}$. Then
\[
W \geq \hat{W}\geq 0, \quad Z \geq \hat{Z}\geq 0,
\]
\[
\hspace{0.5cm}\frac{1}{d_j} \geq \frac{1}{\hat{d}_j} >0,\quad j=1,2,\ldots,n,
\]
where $\tilde{D}={\rm{diag}}(\tilde{d}_1,\ldots,\tilde{d}_n)$. }\\

\noindent {\bf Proof.} The proof of this theorem is quite similar to Proposition 2.2 in \cite{Lee2} and is omitted here. \qquad $\Box$

\bigskip

This theorem shows that the FFAPINV algorithm is well-defined for M-Matrices. We mention that this is correct for the BFAPINV algorithm, as well (see \cite{Lee2}).\\

\noindent \textbf{ Theorem 2.} {\it Let  $A$ be an H-matrix and $\hat{A}$ be its comparison matrix. Let also $A^{-1}=ZDW$ and $\hat{A}^{-1}=\hat{Z}\hat{D}\hat{W}$ be the computed inverse in the factored form by FFINV algorithm for $A$ and $\hat{A}$, respectively. Then
\[
    \vert \frac{1}{d_i } \vert \geq  \frac{1}{\hat{d}_i }>  0.
\]}
\noindent {\bf Proof.} The elements $w_{jl}$ and $z_{lj}$ may be assumed as rational functions
\begin{eqnarray*}
   w_{jl}&=& F_{jl}(a_{11},\ldots,a_{j-1,n},d_{1},\ldots,d_{j-1}),\\
   z_{lj}&=& G_{lj}(a_{11},\ldots,a_{n,j-1},d_{1},\ldots,d_{j-1}).
\end{eqnarray*}
In fact, $w_{jl}$ $(z_{lj})$ is a rational function of first $j-1$ columns ($j-1$ rows) of $A$ and first $j-1$ diagonal entries of $D$. In the same way, $\hat{w}_{jl}$ $(\hat{z}_{lj})$ is a rational function $F_{jl}$ $(G_{lj})$ of first $j-1$ columns ($j-1$ rows) of $\hat{A}$ and first $j-1$ diagonal entries of $\hat{D}$. Let us also assume that
\begin{eqnarray*}
  \tilde{ w}_{jl}=F_{jl}( {\hat a}_{11},\ldots, {\hat a}_{j-1,n},\vert d_{1}\vert,\ldots,\vert d_{j-1}\vert),\\
  \tilde{ z}_{lj}=G_{lj}( {\hat a}_{11},\ldots,{\hat a}_{n,j-1},\vert d_{1}\vert,\ldots,\vert d_{j-1}\vert).
\end{eqnarray*}
This means that $\tilde{ w}_{jl}$ and $\tilde{ z}_{lj}$ are computed similar to  ${\hat w}_{jl}$ and ${\hat z}_{lj}$, with
$\vert d_{1}\vert,\ldots,\vert d_{j-1}\vert$ instead of $d_{1},\ldots,d_{j-1}$.
By induction, we prove that
\begin{equation}\label{e:eq06}
\left\{\begin{array}{l}
a)\quad \vert \frac{1}{d_{k}} \vert \ge \frac{1}{\hat{d}_{k}},\\
b)\quad \hat{w}_{kt} \ge \tilde{w}_{kt} \ge 0, \\
c)\quad \hat{l}_{t}^{(k)}\le \tilde{l}_{t}^{(k)}\le 0,\\
d)\quad \hat{z}_{tk}\ge \tilde{z}_{tk}\geq 0,\\
e)\quad \hat{u}_{t}^{(k)}\le \tilde{u}_{t}^{(k)}\le 0,\\
\end{array}\right.
\end{equation}
for $k=1,\ldots,n$ and $t\le k-1$. Note that $\hat{l}_t^{(k)}$, $\tilde{l}_t^{(k)}$, $\tilde{u}_t^{(k)}$ and $\hat{u}_t^{(k)}$ are defined similar to $\hat{w}_{jl}$ and $\tilde{w}_{jl}$. For $k=1$, there is nothing to prove. Now, let all of these relations hold for every $k \le j-1$. We show that all of them  are correct for $k=j$, as well. From step 4 of Algorithm 2, for every $t\le j-1$, we have
\begin{eqnarray*}
 \hat{l}_{t}^{(j)}\hspace{-0.25cm}&=&\hspace{-0.25cm}\hat{a}_{jt}+\sum^{t-1}_{i=1}\hat{a}_{ji}\hat{z}_{it}.
\end{eqnarray*}
From the hypothesis for every $t\le j-1$, we have $\hat{z}_{it}\ge \tilde{z}_{it}\ge 0$. Therefore
\begin{equation}\label{e:eq07}
  \begin{array}{ll}
    \hat{l}_{t}^{(j)} = \hat{a}_{jt}+\sum^{t-1}_{i=1}\hat{a}_{ji}\hat{z}_{it}\le \hat{a}_{jt}+\sum^{t-1}_{i=1}\hat{a}_{ji}\tilde{z}_{it}=\tilde{l}_{t}^{(j)}\le0.\\
  \end{array}
\end{equation}
Also, from steps 7-9 of Algorithm 2, we have
\begin{eqnarray*}
\hat{w}_{jt}=-\hat{l}_{t}^{(j)}\hat{d}_{t}-\sum^{j-1}_{i=t+1}\hat{l}_{i}^{(j)}\hat{d}_{i}\hat{w}_{it},\qquad t=1,\ldots,j-1.
\end{eqnarray*}
From (\ref{e:eq07}) and  the hypothesis, for every $k\le j-1$, we have
\begin{eqnarray*}
\hat{d}_{k}\ge \vert{d}_{k}\vert\ge 0  &\Rightarrow& -\hat{l}_{k}^{(j)} \hat{d}_{k}\ge -\tilde{l}_{k}^{(j)} \vert{d}_{k}\vert,\\
\hat{w}_{kt}\ge \tilde{w}_{kt}\ge 0 &\Rightarrow &  -\hat{l}_{k}^{(j)} \hat{d}_{k}\hat{w}_{kt}\ge -\tilde{l}_{k}^{(j)} \vert{d}_{k}\vert\tilde{w}_{kt}, \qquad t<k.
\end{eqnarray*}
Hence, we conclude that
\begin{eqnarray*}
  \hat{w}_{jt}=-\hat{l}_{t}^{(j)}\hat{d}_{t}-\sum^{j-1}_{i=t+1}\hat{l}_{i}^{(j)}\hat{d}_{i}\hat{w}_{it}\ge -\tilde{l}_{t}^{(j)}\vert{d}_{t}\vert-\sum^{j-1}_{i=t+1}\tilde{l}_{i}^{(j)}\vert{d}_{i}\vert\tilde{w}_{it}=\tilde{w}_{jt}\ge 0.
\end{eqnarray*}
In the same way, it can be verified that
\begin{eqnarray*}
\hat{u}_{t}^{(j)}\hspace{-0.25cm}&=&\hspace{-0.25cm}\hat{a}_{tj}+\sum^{t-1}_{i=1}\hat{a}_{ij}\hat{w}_{ti}\le \hat{a}_{tj}+\sum^{t-1}_{i=1}\hat{a}_{ij}\tilde{w}_{tj}=\tilde{u}_{t}^{(j)}\le0,\\
  \hat{z}_{tj} \hspace{-0.25cm}&=&\hspace{-0.25cm} -\hat{u}_{t}^{(j)}\hat{d}_{t}-\sum^{j-1}_{i=t+1}\hat{u}_{i}^{(j)}\hat{d}_{i}\hat{z}_{ti}\ge -\tilde{u}_{t}^{(j)}\vert{d}_{t}\vert-\sum^{j-1}_{i=t+1}\tilde{u}_{i}^{(j)}\vert{d}_{i}\vert\tilde{z}_{ti}=\tilde{z}_{jt}\ge 0.
\end{eqnarray*}
Now, we consider two cases $a_{jj}>0$ and $a_{jj}<0$. If $a_{jj}>0$, then from $\hat{w}_{jt}\ge \tilde{w}_{jt}\ge 0$ and $\hat{a}_{tj}\le 0$, we conclude that
\begin{eqnarray*}
 \frac{1}{ \hat{d}_{j}}=\hat{a}_{jj}+\sum^{j-1}_{t=1}\hat{a}_{tj}\hat{w}_{jt}\le \hat{a}_{jj}+\sum^{j-1}_{t=1}\hat{a}_{tj}\tilde{w}_{jt}.
\end{eqnarray*}
Consider the polynomials of $w_{jt}$ and $\tilde{w}_{jt}$. It is easy to see that the corresponding terms of these polynomials have the same absolute values. In other words, they may differ only by the sign. On the other hand, every term of $\sum^{j-1}_{t=1}\hat{a}_{tj} \tilde{w}_{jt}$, when considered as a polynomial in elements of $\hat{A}$, is nonpositive, since all terms of $\tilde{w}$ are nonnegative. Therefore, it is less than or equal to
$\sum^{j-1}_{t=1}a_{tj}w_{jt}$. Since, it is enough that one of its terms to be nonnegative. Putting these results together indicates that
\begin{equation}\label{e:eq10}
 \frac{1}{ \hat{d}_{j}} \le \hat{a}_{jj}+\sum^{j-1}_{t=1}\hat{a}_{tj}\tilde{w}_{jt}\le a_{jj}+\sum^{j-1}_{t=1}a_{tj}w_{jt}=\frac{1}{d_{j}}.
\end{equation}
If $a_{jj}<0$, then from $\hat{w}_{jt}\ge \tilde{w}_{jt}\ge 0$ and $\hat{a}_{tj}\le 0$, we have
\begin{eqnarray*}
 -\frac{1}{ \hat{d}_{j}}=-\hat{a}_{jj}-\sum^{j-1}_{t=1}\hat{a}_{tj}\hat{w}_{jt}\ge -\hat{a}_{jj}-\sum^{j-1}_{t=1}\hat{a}_{tj}\tilde{w}_{jt}.
\end{eqnarray*}
Similar to the case $a_{jj}>0$, one can conclude that
\begin{equation}\label{e:eq11}
 -\frac{1}{ \hat{d}_{j}} \ge  -\hat{a}_{jj}-\sum^{j-1}_{t=1}\hat{a}_{tj}\tilde{w}_{jt}\ge -\hat{a}_{jj}+\sum^{j-1}_{t=1}a_{tj}w_{jt}=a_{jj}+\sum^{j-1}_{t=1}a_{tj}w_{jt}=\frac{1}{d_{j}}.
\end{equation}
Now from Eqs. (\ref{e:eq10}) and (\ref{e:eq11}), we have
\begin{eqnarray*}
 \vert \frac{1}{d_{j}}\vert\ge \frac{1}{ \hat{d}_{j}}.
\end{eqnarray*}
This together with Theorem 1 give the desired result. \qquad $\Box$\\
\\
\noindent \textbf{ Theorem 3.} {\it Let  $A$ be an H-matrix and $\hat{A}$ be its comparison matrix. Let also $A^{-1}\approx ZDW$ and $\hat{A}^{-1} \approx \hat{Z}\hat{D}\hat{W}$ be the factored approximate inverse computed by FFAPINV algorithm for $A$ and $\hat{A}$, respectively. Then
\[
\vert \frac{1}{d_i } \vert \geq  \frac{1}{\hat{d}_i }>  0.
\]}
\noindent {\bf Proof.} The proof is quite similar to the proof of the previous theorem. Indeed, the proof can be done by induction and noting that all of the inequalities have been obtained by comparing  both sides of the inequalities term-by-term. \qquad $\Box$

\bigskip

From this theorem, the following corollary can be easily concluded.

\bigskip

\noindent {\bf Corollary 1.} {\it Let  $A$ be an H-matrix and $D={\rm diag}(d_1,\ldots,d_n)$  be the diagonal matrix computed by Algorithm 1. Then, the sign of $a_{jj}$ and $d_j$ are the same.}

\bigskip

\bigskip

\medskip

\noindent  \textbf{3. An ILU preconditioner based on the FFAPINV algorithm}

\bigskip
\medskip

We first state and prove the following theorem.

\bigskip

\noindent \textbf{Theorem 1.} \textit{Let $W$ and $Z$ be the inverse factors of $A$ in Eq. (\ref{e:eq02}). Also suppose that $L:=W^{-1}$ and $U:=Z^{-1}$. Therefore $A=LD^{-1}U$ is the LDU factorization of $A$ and for $i\le j$
\[
L_{ji}=\beta_{i}^{(j)}, \qquad U_{ij}=\alpha_{i}^{(j)},
\]
in which $\alpha_i^{(j)}$ and $\beta_i^{(j)}$ are computed in steps 5 and 13 of Algorithm 2.}\\

\noindent{\textbf{Proof.}} Let $WAZ=D^{-1}$. Then $A=W^{-1}D^{-1}Z^{-1}$ is the LDU factorization of $A$.
From $Z^{-1}=U$ we conclude $ZU=I$, where $I$ is the identity matrix. By equating the $j$th column of both sides of $ZU=I$, we deduce that
\begin{eqnarray}\label{eq:comput zj}
z_j=e_j-\sum_{i=1}^{j-1} U_{ij}z_i,
\end{eqnarray}
where $e_j$ is the $j$th column of $I$. Since $z_i,i=1,2,\ldots,j-1$, and $e_j$ are linearly independent vectors, then relation (\ref{eq:comput zj}) together with (\ref{e:eq04}) results in $\alpha_{i}^{(j)}=U_{ij}$. In the similar way we can prove that $L_{ji}=\beta_{i}^{(j)}$.  \qquad $\Box$

\bigskip

By the above discussion we propose the next algorithm that computes an ILU factorization of $A$ as by-product of FFAPINV process. We term this algorithm ILUFF (refer to an ILU preconditioning based on FFAPINV algorithm). The algorithm is as follows:

\bigskip

\underline{{\bf Algorithm 3: ILUFF algorithm}}
\vspace{-0.1cm}

\begin{enumerate}
\item  Set $L=U=I_{n \times n}$, $z_1:=e_1$, $w_1:=e_1^T$ and $d_1:=\frac{1}{a_{11}}$ \vspace{-0.25cm}

\item For $j=2,\ldots,n,$~ Do \vspace{-0.25cm}

\item \qquad  $z_j:=e_j$; \quad $w_j:=e_j^T$ \vspace{-0.25cm}

\item \qquad  For $i=1,\ldots,j-1,$ Do \vspace{-0.25cm}

\item \qquad\qquad  $U_{ij}:=d_i w_iA_{*j}$\vspace{-0.25cm}

\item  \qquad\qquad If $|U_{ij}|> \tau$, then $z_j:=z_j-U_{ij} z_i$ \vspace{-0.25cm}

\item \qquad\qquad Drop entries of $z_j$  whose absolute values are smaller than $\tau$\vspace{-0.25cm}

\item \qquad  EndDo \vspace{-0.25cm}

\item \qquad  For $i=1,\ldots,j-1,$ Do \vspace{-0.25cm}

\item \qquad\qquad  $L_{ji}:=d_iA_{j*}z_i$\vspace{-0.25cm}

\item  \qquad\qquad If $|L_{ji}|> \tau$, then $w_j:=w_j-L_{ji}w_i$  \vspace{-0.25cm}

\item \qquad\qquad Drop entries of $w_j$ whose absolute values are smaller than $\tau$\vspace{-0.25cm}

\item \qquad  EndDo \vspace{-0.25cm}

\item \qquad  $d_j:=\frac{1}{w_jA_{*j}}$ \vspace{-0.25cm}

\item  EndDo \vspace{-0.25cm}

\item  Return $L=(L_{ji})$, $U=(U_{ij})$ and $D={\rm diag}(d_1,d_2,\ldots,d_n)$ ($A \approx LD^{-1}U$)
\end{enumerate}

\bigskip

Some consideration can be given here. Obviously, when the matrix $A$ is symmetric then $W=Z^T$ and the computations will be halved. In the case that the matrix $A$  is symmetric positive definite then we have $$d_j=\frac{1}{z_j^TAz_j}>0,$$ and the algorithm is free from breakdown. This is true for nonsymmetric positive definite matrices, as well \cite{Salkuyeh}. Therefore, if the matrix $A$ is positive definite (symmetric or nonsymmetric) then we can use $$d_j:=\frac{1}{z_j^TAz_j},$$ in step 14 of Algorithm 3.

\bigskip

\medskip

\noindent  \textbf{4. Numerical experiments}

\bigskip

\medskip

In this section, we have used the ILUFF as the right preconditioner to solve the linear system of equations (\ref{e:eq01}) with GMRES(50) method. The ILUFF code is written in Fortran 77. But the GMRES(50) code in Sparskit package \cite{Sparskit} has been used. In the first part of the numerical experiments we used 38 nonsymmetric test matrices. All of these matrices have been taken from  the University of Florida Sparse Matrix Collection
\cite{Davis} and none of them are positive definite. In all the experiments whenever a zero pivot has occurred, then we have replaced the zero by the square root of the machine precision. The machine, we used for the experiments, has one quad Intel(R) CPU and 8GB of RAM memory. The initial guess for the iterative solver was always a zero vector. The stopping criterion used was
\begin{equation}\label{e:stopping criterion}
  \frac{\|r_k\|_2}{\|r_0\|_2}< 10^{-10},
\end{equation}
where $r_k$ is the residual of the unpreconditioned system in the $k$th iterate.
We have considered the exact solution as the vector $e=(1,\cdots,1)^T$ and the vector $b=Ae$. We have used the Multilevel Nested Dissection reordering as a preprocessing \cite{karK98}. To implement ILUFF algorithm, it is clear that matrix $A$ should be accessed row-wise and column-wise. In our implementations, we have used just the Compressed Sparse Row (CSR) format of storage \cite{Saadbook} to traverse $A$ row-wise. For the column-wise traverse of $A$, we have used the linked lists \cite{ILUC}.

\par In Table \ref{table:matrix properties}, properties of test matrices and the convergence results of the GMRES(50) method without preconditioning have been reported. In this table, $n$ and $nnz$ are the dimension and the number of nonzero entries of the matrix, respectively. $Its$ stands for the number of iterations and $Time$ is the iteration time which has been computed by $dtime$ command. The times are in second. In this table, a $+$ symbol means that the stopping criterion has not been satisfied after 10,000 number of iterations.

\par Table \ref{table:ILUFF-1}, includes properties of ILUFF preconditioner and results of right preconditioned GMRES(50). Considering Algorithm 3, $\tau$ is the tolerance parameter to drop entries of $L,U$ and $W,Z$ factors. In this table, the parameter $\tau$ has been set to 0.1 for all the test matrices. Suppose that $nnz(A)$, $nnz(L)$ and $nnz(U)$ be the number of nonzero entries of matrix $A$ and $L$ and $U$ factors of ILUFF preconditioner. In the implementation of ILUFF preconditioner, we have merged the $D$ factor into the $U$ factor. Therefore, parameter $density$  in Table \ref{table:ILUFF-1}, is defined as:
\begin{eqnarray*}\label{eq:density}
density=\frac{nnz(L)+nnz(U)}{nnz(A)}.
\end{eqnarray*}
In this table, $Ptime$ stands for the preconditioning time and $It-Time$ is the iteration time to solve the preconditioned system. $Ttime$ is defined as the sum of $Ptime$ and $It-Time$. In this table, $Its$ is the number of iterations of the GMRES(50) to solve the right preconditioned system. Numerical results presented in Table 2, show that the proposed preconditioner greatly reduces the time and iterations for convergence.

\begin{table}
\centering
\caption{ Test matrices properties together with results
 of GMRES(50) without preconditioning.}\vspace{0.4cm}
\label{table:matrix properties}
\small
\begin{tabular}{|l||r|r|c|c|c|}
\hline
Group/Matrix                        & ~~~~$n$~~~~~    & ~~$nnz$~~~   &   ~~~~Time~~~~ &   ~~~~Its~~~~         \\ \hline \hline
HB/fs\_~183\_~1                     & 183     & 998      &    0.00  &     38             \\ \hline
HB/fs\_~183\_~6                     & 183     & 1000     &    0.01  &     36             \\ \hline
Simon/raefsky1                      & 3242    & 293409   &   2.76   &     3588           \\ \hline
Simon/raefsky2                      & 3242    & 293551   &    3.65  &     4790       \\ \hline
Simon/raefsky5                      & 6316    & 167178   &   0.24   &   306              \\ \hline
Simon/raefsky6                      & 3402    & 130371   &   0.68   &   1385             \\ \hline
Muite/Chebyshev3                    & 4101    & 36879    &   +      &     +           \\ \hline
Oberwolfach/flowmeter5              & 9669    & 67391    &    +     &     +           \\ \hline
Rajat/rajat03                       & 7602    & 32653    &    +     &     +           \\ \hline
HB/sherman3                         & 5005    & 20033    &    +     &     +           \\ \hline
Hamm/memplus                        & 17758   & 99147    &    7.22  &     3878          \\ \hline
FEMLAB/poisson3Da                   & 13514   & 352762   &  0.86    &     342         \\ \hline
Botonakis/FEM\_~3D\_~thermal1       & 17880   & 430740   &   0.88   &   283           \\ \hline
FEMLAB/poisson3Db                   & 13514   & 352762   &  12.98   &   620           \\ \hline
Oberwolfach/chipcool0               & 20082   & 281150   &   +      &     +           \\ \hline
Oberwolfach/chipcool1               & 20082   & 281150   &   +      &     +           \\ \hline
Averous/epb1                        & 14734   & 95053    &   2.08   &     1432        \\ \hline
Averous/epb2                        & 25228   & 175028   &   2.81   &     908         \\ \hline
Wang/wang3                          & 26064   & 177168   &   2.34   &     803            \\ \hline
Wang/wang4                          & 26068   & 177196   &   +      &     +           \\ \hline
IBM\_~Austin/coupled                & 11341   & 97193    &   +      &     +           \\ \hline
Simon/venkat01                      & 62424   & 1717792  &   +      &     +           \\ \hline
Sandia/ASIC\_~100ks                 & 99190   & 578890   &  23.85   &    1887             \\ \hline
Hamm/hcircuit                       & 105676  & 513072   &    +     &     +           \\ \hline
Norris/lung2                        & 109460  & 492564   &   +      &     +           \\ \hline
IBM\_~EDA/dc1                       & 116835  & 766396   &   +      &     +           \\ \hline
IBM\_~EDA/dc2                       & 116835  & 766396   &    +     &     +            \\ \hline
IBM\_~EDA/dc3                       & 116835  & 766396   &   +      &     +           \\ \hline
IBM\_~EDA/trans4                    & 116835  & 749800   &   +      &     +           \\ \hline
IBM\_~EDA/trans5                    & 116835  & 749800   &   +      &     +           \\ \hline
Botonakis/FEM\_~3D\_~thermal2       & 147900  & 3489300  &   18.32  &   652           \\ \hline
QLi/crashbasis                      & 160000  & 1750416  &   10.09  &   437            \\ \hline
FEMLAB/stomach                      & 213360  & 3021648  &   11.64  &     344         \\ \hline
Sandia/ASIC\_~320ks                 & 321671  & 1316085  &  10.39   &     201         \\ \hline
Sandia/ASIC\_~680ks                 & 682712  & 1693767  &    13.72 &     84              \\ \hline
Bourchtein/atmosmodd                & 1270432 & 8814880  &    241.92&     707              \\ \hline
Bourchtein/atmosmodl                & 1489752 &	10319760 &   166.14 &    415               \\ \hline
\end{tabular}
\end{table}

\begin{table} \hspace{0.5cm}
\centering
\caption{{\small Properties of ILUFF preconditioner with $\tau=0.1$ and right preconditioned GMRES(50) method}}
\label{table:ILUFF-1}
{\small
\begin{tabular}{|l||c|c|c|c|r|} \hline
 {\rm Matrix}         &   density& ~~Ptime~~    &     ~~It-Time~~   & ~~Ttime~~ & ~~Its~~     \\ \hline\hline
 fs\_~183\_~1         &  0.55  &   0.81    &    0.00    &     0.81   &  10    \\ \hline
 fs\_~183\_~6         &  0.54  &   0.73    &    0.00    &     0.73   &  10    \\ \hline
 raefsky1             &  0.05  &   0.75    &    1.03    &     1.78   &  1092    \\ \hline
 raefsky2             &  0.08  &   0.75    &    1.45    &     2.20   &  1419    \\ \hline
 raefsky5             &  0.20  &   0.75    &    0.01    &     0.76   &  11    \\ \hline
 raefsky6             &  0.15  &   0.75    &    0.01    &     0.76   &  12    \\ \hline
 Chebyshev3           &  0.33  &   0.75    &    0.09    &     0.84   &   245   \\ \hline
 flowmeter5           &  0.74  &   0.76    &    2.37    &     3.13   &   2106   \\ \hline
 rajat03              &  0.78  &   0.82    &    0.25    &     1.07   &   390   \\ \hline
 sherman3             &  0.83  &   0.73    &    0.82    &     1.55   &   1747   \\ \hline
 memplus              &  0.39  &   0.77    &    0.89    &     1.66   &   376   \\ \hline
 poisson3Da           &  0.18  &   0.75    &    0.55    &     1.30   &   180   \\ \hline
 FEM\_~3D\_~thermal1  &  0.22  &   0.76    &    0.23    &     0.99   &   53   \\ \hline
 poisson3Db           &  0.18  &   1.06    &    9.40    &    10.46   &   395   \\ \hline
 chipcool0            &  0.35  &   0.78    &    1.22    &     2.00   &   321   \\ \hline
 chipcool1            &  0.35  &   0.75    &    1.70    &     2.45   &   446   \\ \hline
 epb1                 &  0.78  &   0.73    &    1.11    &     1.84   &   547   \\ \hline
 epb2                 &  0.57  &   0.75    &    0.83    &     1.58   &   209   \\ \hline
 wang3                &  0.85  &   0.76    &    0.84    &     1.60   &   201   \\ \hline
 wang4                &  0.55  &   0.80    &    1.13    &     1.93   &   279   \\ \hline
 coupled              &  0.48  &   0.81    &    0.21    &     1.02   &   149   \\ \hline
 venkat01             &  0.34  &   1.07    &    1.67    &     2.74   &   90   \\ \hline
 ASIC\_~100ks         &  0.79  &   1.01    &    0.35    &     1.36   &   23   \\ \hline
 hcircuit             &  0.76  &   0.99    &    8.95    &     9.94   &   513   \\ \hline
 lung2                &  1.03  &   1.00    &    5.67    &     6.67   &   304   \\ \hline
 dc1                  &  0.65  &  35.26    &    4.75    &    40.01   &   234   \\ \hline
 dc2                  &  0.64  &  34.06    &    2.97    &    37.03   &   143   \\ \hline
 dc3                  &  0.64  &  33.75    &    8.64    &    42.39   &   451   \\ \hline
 trans4               &  0.62  &  21.95    &    2.36    &    24.31   &   128   \\ \hline
 trans5               &  0.63  &  21.80    &    7.52    &    29.32   &   397   \\ \hline
 FEM\_~3D\_~thermal2  &  0.22  &   1.19    &    1.41    &     2.60   &   41   \\ \hline
 crashbasis           &  0.58  &   1.19    &    1.82    &     3.01   &   59   \\ \hline
 stomach              &  0.22  &   1.26    &    2.79    &     4.05   &   72   \\ \hline
 ASIC\_~320ks         &  0.66  &   1.45    &    4.02    &     5.47   &   71   \\ \hline
 ASIC\_~680ks         &  0.61  &   1.97    &    0.68    &     2.65   &   6   \\ \hline
 Bourchtein/atmosmodd &  0.64   &  4.02    &  200.66    &   204.68   &  503    \\ \hline
 Bourchtein/atmosmodl &  0.84   &  5.14    &  100.32    &   105.46   &  209    \\ \hline
 \end{tabular}
}
\end{table}

For the second part of the numerical experiments we consider the matrix $atmosmodj$. This matrix belongs to the Group $Bourchtein$ of matrix collection \cite{Davis}. Dimension of this matrix is $1,270,432$ and it has $8,814,880$ number of nonzero entries. Consider system (\ref{e:eq01}) when the coefficient matrix is $atmosmodj$. Also suppose that $b=Ae$ in which $e=(1,\cdots,1)^T$ is the exact solution. We term this artificial system as the $atmosmodj$ system. Without preconditioning, GMRES(50) method  for this system converges in $1312$ number of iterations in about $447.66$ seconds. In Table \ref{table:ILUFF-2}, we have reported the results of ILUFF preconditioner and the right preconditioned GMRES(50) method for this system. In this table, $density$, $Ptime$, $It-Time$, $Ttime$ and $Its$ have the same meaning as in Table \ref{table:ILUFF-1}. All the $Ttime$ and $Its$ in this table are less than numbers $447.66$ and $1312$, respectively.

\begin{table} \hspace{0.5cm}
\centering
\caption{{\small Properties of ILUFF preconditioner and right preconditioned GMRES(50) for the matrix atmosmodj}}
\label{table:ILUFF-2}
{\small
\begin{tabular}{|c||c|c|c|c|c|} \hline
 {\rm $\tau$}         &   ~~density~~ & ~~Ptime~~   &     ~~It-Time~~   & ~~Ttime~~ & ~~Its~~     \\ \hline\hline
 $0.1$ &   0.64   &   4.04   &   240.81   &  244.85   &  603    \\ \hline
 $0.05$ &  0.76   &   4.42   &   232.66   &  237.08   &  561   \\ \hline
 $0.01$ &  1.04   &   5.17   &   198.06   &  203.23   &  464   \\ \hline
 \end{tabular}
}
\vspace{2cm}
\end{table}

\begin{figure}[t]
\begin{center}
\includegraphics[height=7.8cm,width=9.8cm]{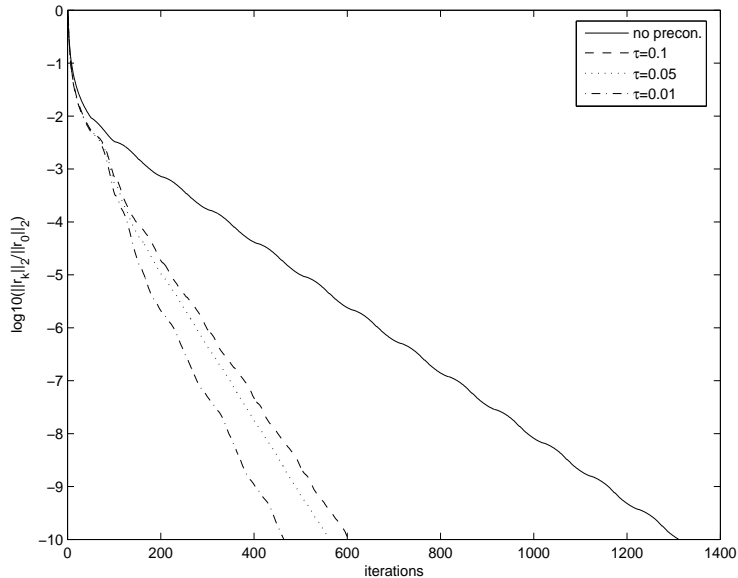}
\end{center}
\caption{\small{\emph{effect of different drop parameters for the ILUFF preconditioner of atmosmodj system.}}}
\label{fig:iteration vs log 2 of related residual}
\end{figure}

 In figure \ref{fig:iteration vs log 2 of related residual}, we have drawn four graphs related to $atmosmodj$ system. In this figure, we take an in-depth look at the results of Table \ref{table:ILUFF-2}.
 The graph with solid line is devoted to the case that GMRES(50) method without preconditioning is applied to solve the $atmosmodj$ system. This graph gives the logarithm of the fraction $\frac{\|r_k\|_2}{\|r_0\|_2}$ for each iterate of the GMRES(50).
 The three other graphs, illustrated by the dashed, dotted and dashed-dotted lines, give the above mentioned logarithm for each iterate $x_k$ of the right preconditioned GMRES(50) method. For these three graphs, the right preconditioner is the ILUFF and has been computed with $\tau$ equal to 0.1, 0.05 and 0.01.
 The shape of the graphs
 first indicates that the ILUFF preconditioner is a good tool to decrease the number of iterations of the GMRES(50) method, applied for $atmosmodj$ system. It also shows the fact that the less the parameter $\tau$ gives the better the quality of the ILUFF preconditioner for $atmosmodj$ system.

\bigskip

\bigskip

\medskip

\noindent  \textbf{5. Conclusion}

\bigskip

\medskip

In this paper, we have proposed an ILU preconditioner based on the Forward FAPINV algorithm say ILUFF which is free from breakdown for nonsymmetric positive definite and H-matrices. Numerical results presented in this paper show that the new preconditioner is very robust and effective.


\bigskip

\bigskip


\begin{thebibliography}{0}

\bibitem{Survey} M. Benzi, {\it Preconditioning techniques for large linear systems: A survey}, J. of Computational
Physics, 182 (2002) 418-477.

\bibitem{Comp} M. Benzi, M. Tuma, {\it A comparative study of sparse approximate inverse preconditioners}, Appl.
Numer. Math., 30 (1999)305-340.

\bibitem{Benzi1} M. Benzi, W. D. Joubert and G. Mateescu,  {\it Numercal experiments with parallel orderings for ILU
preconditioners}, ETNA, 8(1999) 88-114.

\bibitem{Benzi2} M. Benzi, D. B. Szyld and A. van Duin, {\it Ordering for incomplete factorization preconditioning of
nonsymmetric problems}, SIAM J. Sci. Comput., 20 (1999) 1652-1670.

\bibitem{RIF} M. Benzi and M. Tuma, {\it A  robust incomplete factorization preconditioner for positive definite
matrices  approximate inverse preconditioning}, Numer. Linear Algebra Appl., 10 (2003) 385-400.

\bibitem{SAINV}M. Benzi, J. K. Cullum, and M. Tuma, and C. D. Meyer, {\it Robust approximate inverse preconditioning
for the conjugate gradient Method,}  SIAM J. Sci. Comput., 22 (2000) 1318-1332.

\bibitem{SYMAINV}M. Benzi, C. D. Meyer, and  M. Tuma, {\it A sparse approximate inverse preconditioner for the
conjugate gradient method},  SIAM J. Sci. Comput., 17 (1996) 1135-1149.

\bibitem{AINV} M. Benzi, M. Tuma, {\it A sparse approximate inverse preconditioner for nonsymmetric linear systems},
SIAM J. Sci. Comput., 19 (1998) 968-994.

\bibitem{Davis} T. Davis, {\it University of Florida sparse matrix collection}, NA Digest, 92(1994),
\texttt{http://www.cise.ufl.edu/research/sparse/matrices}.

\bibitem{AINVA} S. A. Kharchenko, L. Yu. Kolotilina, A. A. Nikishin, A. Yu. Yeremin, {\it A robust AINV-type method
for constructing sparse approximate inverse preconditioners in factored form}, Numer. Linear Algebra With Appl., 8 (2001) 165–179.

\bibitem{Kol1}L. Y. Kolotilina and A. Y. Yeremin, {\it Factorized sparse approximate inverse preconditioning I.
Theory}, SIAM J. Matrix Anal. Appl., 14 (1993) 45-58.

\bibitem{Kol2}L. Y. Kolotilina and A. Y. Yeremin, {\it Factorized sparse approximate inverse preconditioning II:
Solution of 3D FE systems on massively parallel computers}, Int. J. High Speed Comput., 7 (1995) 191-215.

\bibitem{Lee1} E.-J. Lee and J. Zhang, {\it A two-phase preconditioning strategy of sparse approximate inverse for
indefinite matrices}, Tehnial Report No. 476-07, Department of Computer Science, University of Kentuky, Lexington, KY, 2007.

\bibitem{Lee2} E.-J. Lee and J. Zhang, {\it Fatored approximate inverse preonditioners with dynamic sparsity
patterns}, Tehnial Report No. 488-07, Department of Computer Science, University of Kentuky, Lexington, KY, 2007.

\bibitem{Luo1} J.-G. Luo, {\it An incomplete inverse as a preconditioner for the conjugate gradient method}, Comput.
Math. Appl., 25 (1993) 73–79.

\bibitem{Luo2} J.-G. Luo, {\it A new class of decomposition for inverting asymmetric and indefinite matrices}, Comput.
Math. Appl., 25 (1993) 95–104.

\bibitem{Luo3} J.-G. Luo, {\it A new class of decomposition for symmetric systems}, Mechanics Research Communications,
19 (1992) 159–166.

\bibitem{karK98} G. Karypis and V. Kumar, {\it Fast and High Quality Multilevel Scheme for Partitioning Irregular Graphs}, SIAM J. Sci. Comput., 20(1) (1998) 359-392.

\bibitem{ILUC} N. Li, Y. Saad, E. Chow, \emph{Crout version of ILU for general sparse matrices}, SIAM J. Sci. Comput., 25(2) (2003) 716-728.

\bibitem{Meurant} G. Meurant, {\it The block preconditioned conjugate gradient Method on vector computers,} BIT,
24 (1984) 623-633.

\bibitem{Meij} J. A. Meijerink and H. A. van der Vorst, {\it An iterative solution method for linear systems of which
the coefficient matrix is a symmetric M-matrix}, Math. Comput., 31 (1977) 148-162.

\bibitem{Mant} T. Manteuffel, {\it An incomplete factorization technique for positive definite linear systems},
Math. Comput., 34 (1980) 473-497.

\bibitem{Rezghi} M. Rezghi and S. M. Hosseini, {\it An ILU preconditioner for nonsymmetric positive definite matrices
by using the conjugate Gram-Schmidt process}, Journal of Computational and Applied Mathematics, 188 (2006) 150-164.

\bibitem{GMRES} Y. Saad and M. H. Schultz, {\it GMRES: A generalized minimal
residual algorithm for nonsymmetric linear systems}, SIAM J. Sci.
Statist. Comput., {\bf 7}(1986) 856-869.

\bibitem{ILUT} Y. Saad, {\it ILUT: A dual theshold incomplete LU preconditioner}, Numer. Linear. Algebra Appl.,
1(1994) 387-402.

\bibitem{ILUM} Y. Saad, {\it ILUM: A multi-elimination ILU preconditioner for general sparse matrices}, SIAM J. Sci.
Comput., 174 (1996) 830-847.

\bibitem{BILUM} Y. Saad and J. Zhang, {\it BILUM: Block versions of multi-elimination and multilevel ILU
preconditioner for general linear sparse systems}, SIAM J. Sci. Comput., 20 (1999) 2103-2121.

\bibitem{BILUTM} Y. Saad and J. Zhang, {\it BILUTM: a domain-based multilevel block ILUT preconditioner for
general sparse matrices}, SIAM J. Matrix Anal. Appl., 21 (1999) 279-299.

\bibitem{Saadbook} Y. Saad, {\it Iterative Methods for Sparse linear
Systems}, PWS press, New York, 1995.

\bibitem{Sparskit} Y. Saad, {\it Sparskit and sparse examples}. NA Digest (1994).
\texttt{http://www-users.cs.umn.edu/~saad/software}. Accessed 2010.

\bibitem{Salkuyeh} D.K. Salkuyeh, {\it A Sparse Approximate Inverse Preconditioner for Nonsymmetric Positive
Definite Matrices}, Journal of Applied Mathematics and Informatics, 28 (2010) 1131-1141.

\bibitem{BiCGSTAB} H. A. van der Vorst, {\it Bi-CGSTAB: a fast and smoothly
converging variant of Bi-CG for the solution of nonsymmetric linear
systems}, SIAM J. Sci. Statist. Comput., {\bf 12} (1992) 631-644.


\bibitem{ZhangJi} J. Zhang, {\it A procedure for computing factored approximate inverse}, M.Sc. dissertation,
Department of Computer Science, University of Kentucky, 1999.

\bibitem{Zhang} J. Zhang, {\it A sparse approximate inverse technique for
parallel preconditioning of general sparse matrices}, Appl. Math.
Comput., {\bf 130} (2002) 63-85.

\end{thebibliography}
\end{document}